\documentclass[10pt,amssymb,amsfonts,psfig]{amsart}
\usepackage{epsfig,amsfonts,amssymb,color,amsmath,amscd,}

\newtheorem*{cor1}{Corollary}

\theoremstyle{remark}

\newtheorem*{rem1}{Remark 1}
\newtheorem*{rem2}{Remark 2}

\newtheorem*{example1}{Example 1}
\newtheorem*{example2}{Example 2}

\newtheorem*{conjecture1}{Conjecture 1} 

\theoremstyle{definition}

\newtheorem*{thm1}{Theorem 1}
\newtheorem*{thm2}{Theorem 2}
\newtheorem*{thm3}{Theorem 3}
\newtheorem*{thm4}{Theorem 4}
\newtheorem*{thm5}{Theorem 5}

\newtheorem*{LY Theorem}{Lee-Yang Theorem}
\newtheorem*{thml}{Theorem (Lyapunov)}
\newtheorem*{lem1}{Lemma 1}
\numberwithin{equation}{section}
\numberwithin{figure}{section}

\font\nt=cmr7

\def\note#1
{\marginpar
{\nt $\leftarrow$
\par
\hfuzz=20pt \hbadness=9000 \hyphenpenalty=-100 \exhyphenpenalty=-100
\pretolerance=-1 \tolerance=9999 \doublehyphendemerits=-100000
\finalhyphendemerits=-100000 \baselineskip=6pt
#1}\hfuzz=1pt}

\renewcommand{\mod}{\operatorname{mod}}

\renewcommand{\Re}{\operatorname{Re}}

\newcommand{\eps}{{\varepsilon}}

\newcommand{\R}{{\Bbb R}}

\def\B0{{\mathbf{0}}}

\catcode`\@=12

\def\Empty{}
\newcommand\oplabel[1]{
  \def\OpArg{#1} \ifx \OpArg\Empty {} \else
  	\label{#1}
  \fi}
		
\newcommand{\comm}[1]{}
\newcommand{\comment}[1]{}

\begin{document}

\bigskip\bigskip

\title[On zeros of the Alexander polynomial of an alternating knot]{On zeros of the Alexander polynomial of an alternating knot}

\author {Lilya Lyubich  and Kunio Murasugi}
\date{\today}

 \begin{abstract}  
We prove that for any zero $\alpha $ of the Alexander 
polynomial of a two-bridge knot, $-3 < \Re(\alpha ) < 6 $. Furthermore, for a large
class of two-bridge knots we prove $-1<\Re(\alpha )$.
\end{abstract}

\setcounter{tocdepth}{1}
 
\maketitle
\tableofcontents

\section{Introduction}

In 2002 Jim Hoste made the following conjecture based on his extensive
computer experiment:
\begin{conjecture1} ( J. Hoste, 2002)
Let $ K $ be an alternating knot and $ \Delta _K (t) $ be its  Alexander 
polynomial. Let $\alpha $ be a zero of $ \Delta _K (t) $. 
Then $\Re(\alpha ) > -1.$
\end{conjecture1}
This conjecture is known to be true for some classes of alternating knots.
\medskip

\noindent 
1) If $K$ is a special alternating knot, then all zeros of its Alexander polynomial lie on a unit circle ([M2],[L],[T]), and 
$\Delta_K(-1) \neq 0,$ so  Conjecture 1 holds.
\medskip

\noindent
2) If $ \alpha $ is a real zero of the Alexander polynomial $ \Delta _K (t) $ of an alternating knot $K$, 
then $\alpha > 0 $, since the coefficients of the Alexander polynomial 
of an alternating knot
have alternating signs ([C],[M1]).
Therefore, if all zeros are real, then K satisfies Conjecture 1.
\medskip

\noindent
3) Any knot $K$ with $ \deg \Delta _K (t)=2 $ satisfies 
$ -1<\Re(\alpha )<3.$ Any alternating knot $K$ with $\deg \Delta _K(t) =4 $
satisfies Conjecture 1.
\medskip

The problem of finding a lower or upper bound of the real part of zeros
 of  the Alexander polynomial is reduced to a problem of showing
 the stability  of the matrix associated to a Seifert matrix $U$ of a knot. Then we apply a well known Lyapunov theorem on the stability of 
matrices. This approach, described in detail in section~2 below, is particularly successful for two-bridge knots. A two-bridge knot $K=K(r)$
is identified by a rational number $r$. We use an even negative 
continued fraction expansion $r=[2a_1,2a_2,\ldots,2a_m] $ to construct
a knot diagram $\Gamma (K(r)) $,  a Seifert surface $F$ and its Seifert
matrix $U$.

Throughout the paper by a two-bridge knot we will mean a two-bridge knot or a two-component two-bridge link, and its Alexander polynomial is defined by
$\Delta _{K(r)}=\det(Ut-U^T)$ (see \cite{BZ}).

In this paper we prove the following theorems:
\begin{thm1} Let $K(r)$ be a two-bridge knot, $\Delta_K(t)$ be its Alexander polynomial and $\alpha $ be a zero of $\Delta _{K}(t). $ Then
\begin{equation*}
-3 < \Re(\alpha ) < 6.
\end{equation*}
\end{thm1}
\begin{thm2} Let $K(r)$ be a two-bridge knot, 
$r=[2a_1,2a_2,\ldots,2a_m]. $ If $a_ia_{i+1} < 0 $ for $i=1,2,\ldots,m\!-\!1$, 
then all zeros are real, hence the conjecture holds.
\end{thm2}
\begin{thm3}Let $K(r)$ be a two-bridge knot, 
$r=[2a_1,2a_2,\ldots,2a_m]. $ If among $a_1,\ldots,a_m $ there are no two
consecutive $1$ or $-1 \;\;$(namely, $a_ia_{i+1}\neq 1 \;$
 for $i=1,\ldots,m\!-\!1)$,
then the conjecture holds.
If moreover $|a_i|>1 $ for $i=1,2,\ldots,m$, then $ -1<\Re{\alpha }<3$.
\end{thm3}

It is known that $K(r) $ is fibered if and only if $|a_j|= 1 
\text{ for all } j.$ 
\begin{thm4} Let $K(r)$ be a fibered two-bridge knot with
\begin{equation*}
\begin{array}{ccccc}
r=&[\;\underbrace{ 2, \ldots , 2,} &  \underbrace{-2, \ldots, -2},&
\ldots,&
\underbrace{(-1)^{m-1}2,\ldots , (-1)^{m-1}2}\;] \\
&k_1 &  k_2&& k_m
\end{array}.
\end{equation*}
If $k_j=1$ or $2$ for all $j$, then the conjecture holds.
\end{thm4}
\begin{thm5}
Let $K(r) $ be a two-bridge knot, $r=r(m,c)=[2c,-2c,\ldots ,(-1)^{m-1}2c]$, $c>0,\;m\geq 1.\;$
Then all zeros of $\;\Delta _{K(r)}\; $ satisfy inequality:
\begin{equation*}
(\dfrac{\sqrt{1+c^2}-1}{c})^2<\alpha <(\dfrac{\sqrt{1+c^2}+1}{c})^2.
\end{equation*}
\end{thm5}
 For non-alternating knots there are no such bounds.
\begin{example1}
Let $\;\Delta _K(t)=1+at-(2a+1)t^2+at^3+t^4,\;\;a>0 .$
Since $\Delta _K(-(a+1))<0$,
there is a zero $\;\alpha \;$ of $\;\Delta _K(t)\;$ such that
$\;\Re (\alpha )<-a-1.$ $K$ is not alternating.
\end{example1}
\begin{example2} Let $\Delta _K(t)=1-2at+(4a-1)t^2 -2at^3+t^4,\;\;
a \geq 4$. Then $\Delta _K(a)<0 $ and hence, there exists a zero 
$\alpha $ such that $\alpha >a $.
 $K$ is not alternating. In fact, if K is alternating, then K is fibered and
since $\deg \Delta _K(t) = 4,\; K\;$  has at most 8 crossings. 
However, such an
alternating knot (including non-prime alternating knots) does not exist in
the table  if $ \;a \geq 4\;$ (see \cite{BZ}).
\end{example2}

\section{ Stability of matrices and Lyapunov theorem } 
Let $K$ be an alternating knot (or link) and 
$ \Delta _K (t)=c_0 + c_1t+ c_2t+ \ldots +c_nt^n $,  $c_n\!\neq \!0$
be its Alexander polynomial. 
Let $A $ be a companion matrix of $\Delta _K (t)$ i.e.
$\Delta _K (t) = c_n \det (tE-A). \; $ The eigenvalues of $A$ are the zeros of $\Delta _K (t).$ We have
\begin{equation*}
\Re(\alpha)>-1 
\Longleftrightarrow \Re(-(1+\alpha )) <0.
\end{equation*}
\noindent Let $\alpha _1, \alpha _2, \ldots ,\alpha _n $  be  all zeros
of $\Delta _K(t)$ (= all eigenvalues of $A).$ 
 Then it is easy to see that $-(1+\alpha _1), -(1+\alpha_2 ),
\ldots ,-(1+\alpha _n), $ are eigenvalues of $-(E + A).$
To prove that all eigenvalues of a matrix have  negative real parts, we apply 
 the Lyapunov theorem: 

\noindent
Let $ M$ be a real $n \times n $ matrix.
Consider a linear vector differential equation
\begin{equation*}
\mathbf{\dot{x}}=M \mathbf{x} .
\end{equation*}  
It is a known theorem in ODE that
all solutions  $ \bold{x}(t) \in \R^n $ of it  are stable, namely 
$\bold{x}(t) \longrightarrow 0 $ as 
$ t \longrightarrow \infty $,
if and only if all eigenvalues of $M$ have negative real parts.
In this case $M$ is called stable.

\begin{thml} \cite{Gant}
All eigenvalues of $ M $ have negative real parts
if and only if there exists a symmetric positive definite matrix $ V $
such that
\begin{equation*}\label{posdef}
VM+M^TV=-W,\; \text{ where } W  \text{ is positive definite}.
\end{equation*}
\end{thml}

\noindent Hence
$K$  satisfies Conjecture 1 if there exists a positive definite
matrix $ V $ such that 
\begin{equation}\label{posdef-1}
V(E+A)+(E+A^T)V  = W \text{ is positive definite.}
\end{equation}

Similarly to (\ref{posdef-1}),
all zeros of $\Delta _K(t) $ satisfy $\; -k< \Re(\alpha ) $ if and only if
$\; -(kE+A) $ is stable, i.e 
there exists a positive definite
matrix $ V $ such that 
\begin{equation*}\label{posdef-k}
V(kE+A)+(kE+A^T)V = W \text{ is positive definite.}
\end{equation*} 

Further, all zeros of $\Delta _K(t) $ satisfy  $\Re (\alpha )<q $ if and only if
$A-qE$ is stable, i.e.
there exists a positive definite
matrix $ V $ such that 
\begin{equation*}\label{posdefq}
V(qE-A)+(qE-A^T)V= W \text{ is positive definite.}
\end{equation*} 
 
To prove that a matrix is positive definite we use 
the following lemma. 
\begin{lem1}(Positivity Lemma)
\begin{equation*}
\text{Let } N= \left [
\begin{array}{lllcc} 
a_{11}&a_{12}&&&\\
a_{21}&a_{22}&a_{23}&&\\
&\ddots & \ddots &\ddots&\\
&&\ddots &\ddots&a_{n-1,n}\\
&&&a_{n,n-1}&a_{n,n}
\end{array}
\right ] \;\;\text{ be a real symmetric matrix. }
\end{equation*}
Suppose that for $1\leq j \leq n$

\noindent
 (i)$\; a_{j,j} > 0 ,\; \;a_{j,j-1}, a_{j,j+1} \neq 0,\;$
and all non-specified entries are $0$.

\noindent
(ii)$\; a_{j,j} \geq |a_{j,j-1}|+ |a_{j,j+1}|,$

\noindent
(iii) there exists $ i \text{ such that }\; a_{i,i} > |a_{i,i-1}|+ |a_{i,i+1}|.$

\noindent
Then N is positive definite.
\end{lem1}

\noindent
The proof is by induction.

\section{ Two-bridge knots}

Let $K(r), 0<r=\beta /\alpha <1, \;\;0<\beta <\alpha $, be a two-bridge
knot or a (two-component) two-bridge link of type $(\alpha ,\beta )$. 
 We can assume one of $\alpha $ and $\beta $ is even. Consider an even (negative)
 continued fraction
expansion of $r$ :

\noindent

\noindent
\begin{align*}
r= \beta /\alpha = \cfrac{1}{\displaystyle 2a_1 - 
\cfrac{1}{\displaystyle 2a_2- } }- && \\
&\ddots  \\
& \hspace{10mm} -\frac{1}{\displaystyle 2a_m}\\
 =[2a_1, 2a_2,\ldots ,2a_m]  &.\\ 
\end{align*}

\noindent 
This expansion is unique.
We obtain from it
a knot or a link diagram $\Gamma (K(r)) \text{ of } K(r).$ 
(see Fig.1)

\begin{figure*}[h]
\begin{align*}
K(r) \text{ is a knot} \hspace{10mm} K(r) \text{ is a link}\\
m=0(\mod 2) \hspace{10mm} m=1(\mod 2)
\end{align*}

\includegraphics[height=10cm]{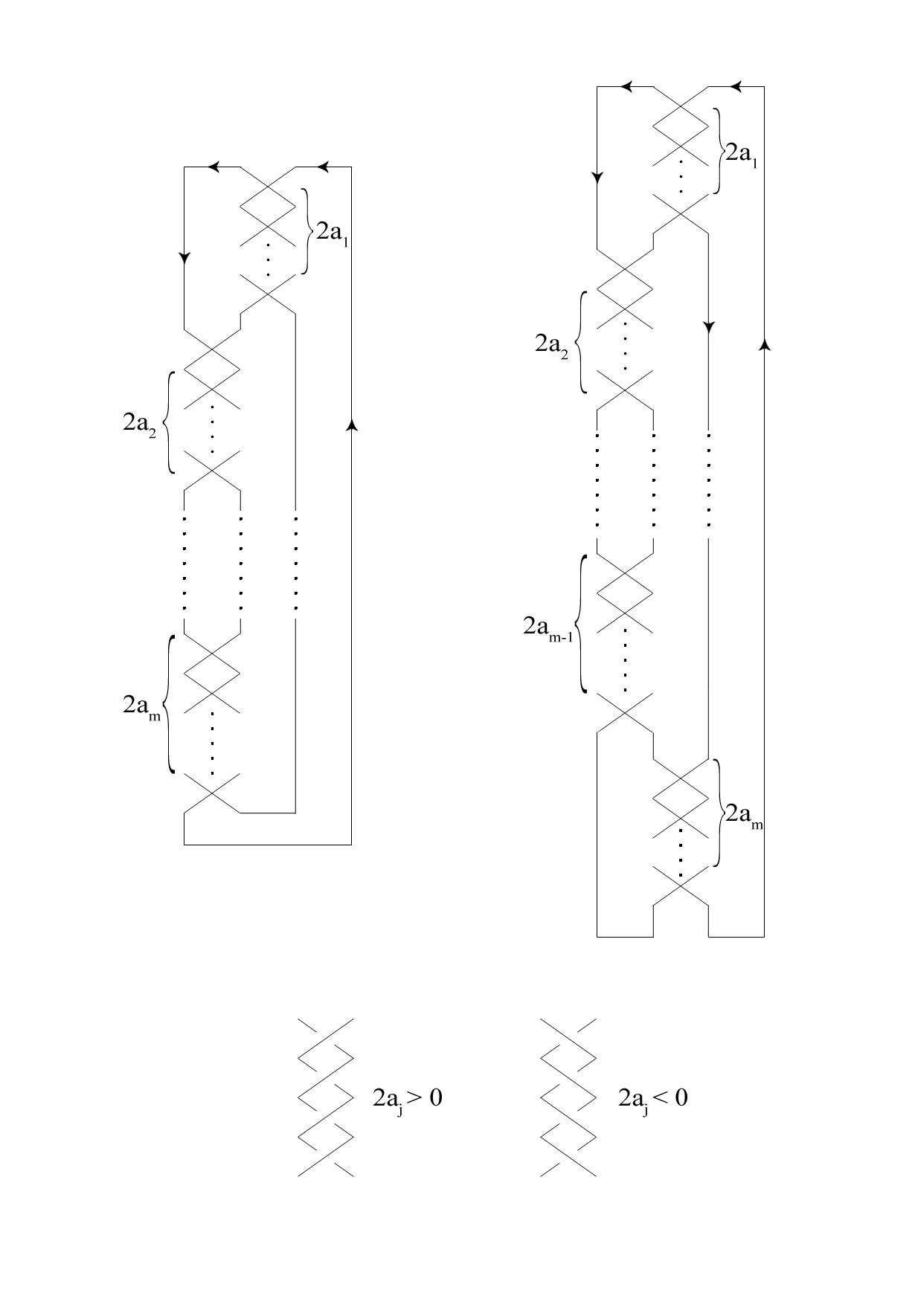}

Figure 1.
\label{knot}
\end{figure*}

The following facts are well known:

\noindent
(1) $ K(r) $ is special alternating if and only if
$a_1, a_2, \ldots , a_m $ are either all positive or all negative.

\noindent
(2)$K(r) $ is fibered if and only if $|a_j|= 1 \text{ for all } j.$  

\noindent 
(3)$\Gamma (K(r))$ is an alternating diagram if and only if 
$ a_j a_{j+1}< 0 \text{ for } j=1,2, \ldots m-1.$ 

\noindent
(4) $\Gamma (K(r)) $ gives a minimal genus Seifert surface $ F $
for $ K(r) \;$(see Fig.2).

\begin{figure}[h]
\includegraphics[height=60mm]{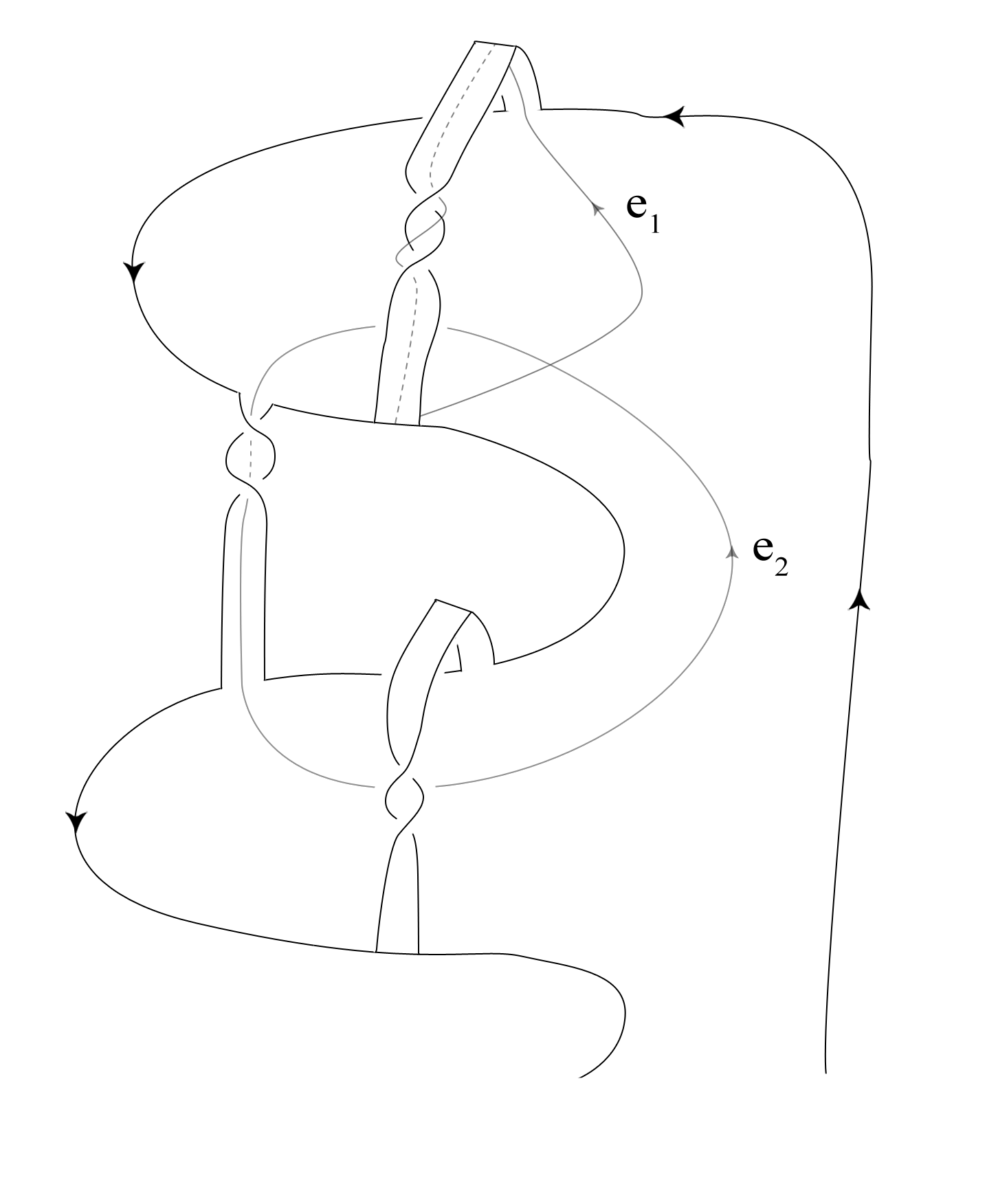}

Figure 2. Seifert surface F.
\label{seifert}
\end{figure}

We use this Seifert surface to calculate a Seifert matrix $\;U=(u_{ij})\;$
of $K,\; u_{ij}=lk(e_i^{\#} ,e_j),\; i,j=1,\ldots ,m. $
For the fragment of $ F $ with only two bands with (half)twists $2a_1 $ 
and $ 2a_2 $   we have 
\begin{equation*}
\begin{array}{cc}
lk(e_1^{\#} ,e_1) = a_1 ,&\;\;lk(e_1^{\#} ,e_2)=0,  \\
\vspace{2mm}
lk(e_2^{\#} ,e_1)  = -1 ,&\;\;lk(e_2^{\#} ,e_2) =a_2, 
\end{array}
\end{equation*}
and in general, it is not difficult to see that for a two-bridge knot 
$K=[2a_1,2a_2,\ldots, 2a_m]$  a Seifert matrix corresponding 
to the surface $ F $ is:
\vspace{2mm}
\begin{equation}\label{U}
U= \left[
\begin{array}{ccccccc}
a_1 & 0 & & & & \\
-1& a_2 & 1& & & \\
& 0 & a_3 & 0 &  & \\
& & -1 & a_4 & 1 & \\
& & & & \ddots & \\
& & & &  -1 & a_m 
\end{array}
\right ]\; \; \text{ or }
\;\; \left [
\begin{array}{ccccccc}
a_1 & 0 & & & & \\
-1& a_2 & 1& & & \\
& 0 & a_3 & 0 &  & \\
& & -1 & a_4 & 1 & \\
& & & & \ddots & \\
& & & & 0 & a_m 
\end{array}
\right ] 
\end{equation}

\noindent
(depending on $\;m\;$ being even or odd, respectively), 
where all non-specified entries are $ 0 $.
The Alexander polynomial of $K$ is $ \Delta_K(t)= \det(tU-U^T).$ 
 So $ A=U^{-1}U^T $ is a companion matrix for $\Delta_K(t).\;$
We have
\begin{equation}
U^{-1}= \left [
\begin{array}{ccccc}
\frac{1}{a_{1}}&&&\ldots 0 & \\
\vspace{2mm}
\frac{1}{a_1a_2}& \frac{1}{a_2} &  -\frac{1}{a_2a_3}&& \\
\vspace{2mm}
&&\frac{1}{a_3} &&\\
\vspace{2mm}
&&\frac{1}{a_3a_4}& \frac{1}{a_4} &  -\frac{1}{a_4a_5}\\
\vspace{2mm}
&0 &&& \ddots
\end{array}
\right ]\hspace{5mm} 
\end{equation}
\vspace{2mm}
and
$U^{-1}U^T=$~ 
\begin{equation}\label{U^-1U^T}
= \left [
\begin{array}{clclcccc}
1\hspace{2mm} & -\frac{1}{a_1} & 0 &&\ldots &&&\\
\vspace{2mm}
\frac{1}{a_2} \hspace{2mm}& 1\!-\!\frac{1}{a_1a_2}\!-\!\frac{1}{a_2a_3}&\hspace{2mm}
-\frac{1}{a_2}\hspace{2mm} & \frac{1}{a_2a_3}& 0 & \ldots & &\\
\vspace{2mm}
0&\frac{1}{a_3}& 1 & -\frac{1}{a_3}& 0 & 0 & 0 & \ldots \\
\vspace{2mm}
 0 & \frac{1}{a_3a_4} & \frac{1}{a_4} & 1\!-\!\frac{1}{a_3a_4}
\!-\!\frac{1}{a_4a_5}&\hspace{2mm} -\frac{1}{a_4} & \frac{1}{a_4a_5}& 0 & \ldots  \\
\vspace{2mm}
0& 0 &0& \frac{1}{a_5}& 1 & -\frac{1}{a_5}& 0& \ldots\\
\vspace{4mm}
&&& \ldots &&&\\
\end{array}
\right ]
\end{equation}

\noindent
The last row of $\;A\;$ is 
$\;[0,\ldots ,0,\frac{1}{a_{m-1}a_m},
\frac{1}{a_m},1\!-\!\frac{1}{a_{m-1}a_m}]\;$ if
$\;m\;$ is even, and
$\;[0,\ldots ,0,\frac{1}{a_m},1]\;$  if $\;m\;$ is odd.

\vspace{2mm}

\noindent 
$\;A=U^{-1}U^T \;$  is a companion matrix for the Alexander polynomial 
of the two-bridge knot
$K(r)$, where 
$r=[2a_1,2a_2, ... 2a_m] $.

\section{Theorem 1: Lower and upper bounds on the real part of zeros for two-bridge knots}

In this section we prove the following theorem: 
\begin{thm1}
If $\alpha $ is a zero of the Alexander polynomial of a two-bridge knot, then
\begin{equation*}
-3<\Re(\alpha )<6.
\end{equation*}
\end{thm1}
\begin{proof}
a) To show that $\Re (\alpha )> -k $ we prove that $-(kE+A) $ is stable.
Taking $V=E$, it is enough to show that 
$A_0=(kE+A) + (kE+A^T)=2kE+A+A^T $ is positive definite.  
Now, $ A_0 $ is of the form
\vspace{2mm} 
\begin{equation*}
A_0= \left [
\begin{array}{cccccccc}
\vspace{2mm}
2k \! +\!2& \; b_1 &&&&&&\\
\vspace{2mm}
b_1 & \; c_1 \; &b_2 \; &d_1 \; &&&&\\
\vspace{2mm}
& b_2 \;& 2k \!+ \! 2 \;& b_3 \;&&&&\\
\vspace{2mm}
& d_1 \;& b_3 \; & c_2 \; & b_4 \; & d_2 &&\\  
\vspace{2mm}
&&&b_4 \;& 2k \!+\!2 \;&b_5\;&&\\
\vspace{2mm}
&&& d_2 \;& b_5 \;& c_3\;& b_6 \;& d_3\\ 
\vspace{2mm}
&&&&& b_6 \;& 2k \!+ \!2 \;&b_7 \\
\vspace{2mm}
&&&&&&&\ddots
\end{array}
\hspace{10mm} \right ]
\end{equation*}
where
$
l=[\frac{m}{2}],$

\noindent
$ b_j=-\frac{\displaystyle 1}{\displaystyle a_j}+
\frac{\displaystyle 1}{\displaystyle a_{j+1}}, \; \;j=1,\ldots ,m-1 ,$

\noindent
$
c_j =(2k+2)-\frac{\displaystyle 2}{\displaystyle a_{2j-1}a_{2j}} - 
\frac{\displaystyle 2}{\displaystyle a_{2j}a_{2j+1}}, \;\;  
j=1,\ldots ,l-1,
$
\vspace{1mm}

\noindent
$
c_l=(2k+2)-\frac{\displaystyle 2}{\displaystyle a_{2l-1}a_{2l}} \;\;$ 
for $ \;m \;$
even,
\vspace{1mm}

\noindent
$
c_l=(2k+2)-\frac{\displaystyle 2}{\displaystyle a_{2l-1}a_{2l}}-
\frac{\displaystyle 2}{\displaystyle a_{2l}a_{2l+1}} \;\;$ for $ m \;$ odd.

\noindent
$
d_j=\frac{\displaystyle 1}{\displaystyle a_{2j}a_{2j+1}} + 
\frac{\displaystyle 1}{\displaystyle a_{2j+1}a_{2j+2}} ,\;\;
j=1,\ldots  ,l-1. 
$
\vspace{1mm}

\begin{equation*}
\text{Let }  P =
\left [
\begin{array}{ccccccc}
1 & & & & & &\\
-\dfrac{ b_1}{ 2k+2} &  1 &
-\dfrac{ b_2}{ 2k+2} &   &  &&\\
&& 1 &&&&\\
& &-\dfrac{ b_3}{ 2k+2} & 1 &
-\dfrac{ b_4}{ 2k+2} &&\\
&&& &1&&\\
&&& &&\ddots &\\
\end{array}
\right]
\end{equation*}

\vspace{4mm}
Then 
\begin{equation*}
PA_0P^T =
\left [
\begin{array}{ccccccc} 
2k+2 & && && & \\
&\alpha _1& 0&\beta _1
& && \\ 
& 0 &2k+2 & 0 & &  & \\
& \beta _1 &  0
 & \alpha _2 &0
& \beta _2&\\
&&&0&2k+2 &0&\\
&&& \beta _2&0&
\ddots& \\
&&& \ldots &&&
\end{array}
\right ]
\end{equation*}
\vspace{2mm}
\begin{equation*}\label{oplus}
\approx
\left [
\begin{array}{cccc}
2k+2 &&&0\\
&2k+2 &&\\
\vspace{2mm}
&& \ddots &\\
0 &&&2k+2
\end{array}
\right ]
\oplus
\left [
\begin{array}{cccccc}
\alpha _1 & \beta _1 &&&0&\\
\beta_1& \alpha _2 & \beta_2 &&  &\\
& \beta_2& \alpha _3 & \beta_3 && \\
&&\ddots &\ddots & \ddots &\\
 0 &&&&&
\end{array}
\right ]
\end{equation*}
( denote the second matrix by $A_{00} $),

\noindent
where 
\begin{equation*}\label{alpha}
\alpha _j=
-\frac{ b_{2j-1}^2}{ 2k+2} +c_j
-\frac{b_{2j}^2}{ 2k+2},\;\;\;j=1,\ldots l-1, 
\end{equation*}
\begin{equation*}
\alpha _l=\begin{cases}

  -\frac{ b_{2l-1}^2}{ 2k+2} +c_l &  m \text{ is even}\\
 -\frac{ b_{2l-1}^2}{ 2k+2} +c_l
-\frac{b_{2l}^2}{ 2k+2}&   m  \text{ is odd }.
\end{cases}
\end{equation*}
\begin{equation*}
\beta _j=d_j -\frac{ b_{2j}b_{2j+1}}{ 2k+2} ,\;\; \; j=1,\ldots ,l-1,\; 
\end{equation*}

\noindent
We show $\;(i)\;\alpha _j>0,\;\; (ii)\; \alpha_j \geq |\beta _{j-1}| +
|\beta_j|\;,\; (iii)\; $ there exists $i$ such that $\alpha _i > |\beta _{i-1}| +
|\beta_i| $.
Then $A_{00} $ is positive definite.

Let $ k=3 .\;$ Then 
$$
 \alpha _j =8-\dfrac{2}{a_{2j-1}a_{2j}}-
\dfrac{2}{a_{2j}a_{2j+1}}-\dfrac{1}{8}\left ( \dfrac{-1}{a_{2j-1}}+
\dfrac{1}{a_{2j}}\right ) ^2
-\dfrac{1}{8}\left ( \dfrac{-1}{a_{2j}}+\dfrac{1}{a_{2j+1}}\right )^2
$$
$$
=8-\dfrac{12}{8} \left( \dfrac{1}{a_{2j-1}a_{2j}}+\dfrac{1}    
{a_{2j}a_{2j+1}}  \right )- 
\dfrac{1}{8}\left ( \dfrac{1}{a_{2j-1}}+\dfrac{1}{a_{2j}}\right ) ^2-
\dfrac{1}{8}\left ( \dfrac{1}{a_{2j}}+\dfrac{1}{a_{2j+1}}\right ) ^2.
$$
\vspace{2mm}
Since $|a_j|\geq 1 $ for all $ j, \;\left | \dfrac{1}{a_j}+
\dfrac{1}{a_{j+1}} \right | \leq 2 \;\text{ and }
\left | \dfrac{1}{a_{2j-1}a_{2j}}+\dfrac{1}{a_{2j}a_{2j+1}}\right |
\leq 2 $
 
\noindent 
and hence  $ \alpha _j \geq 8-\dfrac{3}{2} \cdot 2-\dfrac{1}{8}
\cdot 4 - \dfrac{1}{8} \cdot 4 =4.$

\noindent
On the other hand,
$
\beta _{j-1}=d_{j-1}-\dfrac{b_{2j-2}b_{2j-1}}{8}$

\noindent
$
=\dfrac{1}
{a_{2j-2}a_{2j-1}} +\dfrac{1}{a_{2j-1}a_{2j}}-
\dfrac{1}{8} \left ( -\dfrac{1}{a_{2j-2}}+\dfrac{1}{a_{2j-1}} 
\right )
\left ( -\dfrac{1}{a_{2j-1}}+\dfrac{1}{a_{2j}} 
\right )
$

\noindent
$
=\dfrac{7}{8}\left (\dfrac{1}
{a_{2j-2}a_{2j-1}} +\dfrac{1}{a_{2j-1}a_{2j}} \right )+
\dfrac{1}{8}\left ( \dfrac{1}{a_{2j-2}a_{2j}}+\dfrac{1}{a_{2j-1}^2}
\right ).
$

\noindent
Since $|a_j|\geq 1,\;\; |\beta _{j-1}|\leq \dfrac{7}{8} \cdot 2+
\dfrac{1}{8} \cdot 2 =2 . $ 

\noindent
Similarly $|\beta _j| \leq 2.$ 
Thus $\alpha _j \geq |\beta _{j-1}|+|\beta _j|$
and $\alpha _1 >|\beta _1|.$
If $\beta _j=0 $
then $\alpha _{j+1} >|\beta _{j+1}|.$
This proves the left inequality.

\noindent
b) To prove that $ \Re (\alpha )< q \;$ it is enough to show that
$B_0=(qE-A)+(qE-A^T)=2qE-(A+A^T) \; $ is positive definite. 
$\;B_0 $ is of the form
\begin{equation*}
B_0= \left [
\begin{array}{ccccccccc}
2q-2& -b_1 &&&&&&&\\
-b_1&e_1&-b_2&-d_1&&&&&\\
&-b_2&2q-2 &-b_3 &&&&&\\
&-d_1 &-b_3 &e_2 &-b_4&-d_2&&&\\
&&& -b_4& 2q-2 &-b_5 &&&\\
&&&-d_2 &-b_5&e_3 &-b_6 &-d_3&\\
&&&&&&&\ddots &
\end{array}
\right ]
\end{equation*}

\noindent
where  $e_j =2q-2+\dfrac{2}{a_{2j-1}a_{2j}} +\dfrac{2}{a_{2j}a_{2j+1}},
\;\; j=1, \ldots l-1 ,$

\noindent 
$e_l =2q-2+\dfrac{2}{a_{2l-1}a_{2l}} +\dfrac{2}{a_{2l}a_{2l+1}},\;\; 
\text{ if } m \text{ is odd,}$

\noindent 
$e_l =2q-2+\dfrac{2}{a_{2l-1}a_{2l}} \;\; \text{ if } m \text{ is even, }$

\noindent ($l=[\frac{m}{2}],\;\; $as before.)

\noindent
Using 
\begin{equation*}
Q=\left [
\begin{array}{cccccc}
1&&&&&\\
\dfrac{b_1}{2q-2} &1 & \dfrac{b_2}{2q-2}&&&\\
&& 1 &&&\\
&&\dfrac{b_3}{2q-2}&1&\dfrac{b_4}{2q-2}&\\
&&&&1& \ddots 
\end{array}
\right ]
\end{equation*}
we obtain   
\begin{equation*}
QB_0Q^T=
\left [ 
\begin{array}{ccc}
2q-2&&0\\
&\ddots&\\
0&&2q-2
\end{array}
\right ] 
\oplus \left [
\begin{array}{ccccc}
\gamma _1&\delta _1 &&&\\
\delta _1 &\gamma _2&\delta _2&&\\
&\delta _2 &\gamma _3&\delta _3&\\
&&&\ddots 
\end{array}
\right ],
\end{equation*}
\vspace{2mm}

\noindent
where $\;\gamma _j=e_j-\dfrac{b_{2j-1}^2}{2q-2}-\dfrac{b_{2j}^2}{2q-2},\;\; 
j=1,\ldots l \;\;\text{ for } m \;\text{ odd, }$

\noindent
and  $\;\gamma _l \;$ is replaced by 
$\;\gamma _l= e_l-\dfrac{b_{2l-1}^2}{2q-2}\;\;$ for $\;m \;$ even,

\noindent
$\delta _j=-d_j-\dfrac{b_{2j}b_{2j+1}}{2q-2}, \;\; j=1, \ldots l-1.$

Let $q=6$. Then 

\noindent
$\gamma _j=10+\dfrac{2}{a_{2j-1}a_{2j}}+\dfrac{2}{a_{2j}a_{2j+1}}-
\dfrac{1}{10} \left ( -\dfrac{1}{a_{2j-1}}+\dfrac{1}{a_{2j}} \right ) ^2 - 
\dfrac{1}{10} \left ( -\dfrac{1}{a_{2j}}+\dfrac{1}{a_{2j+1}} \right ) ^2
$
$\geq 10-2-2-\dfrac{1}{10} \cdot 4-\dfrac{1}{10} \cdot 4=5.2,$
since $\left | -\dfrac{1}{a_k} + \dfrac{1}{a_{k+1}} \right | \leq 2.$
While $ \delta _{j-1}=-\dfrac{1}{a_{2j-2}a_{2j-1}}-\dfrac{1}{a_{2j-1}a_{2j}}-
\dfrac{1}{10}\left ( -\dfrac{1}{a_{2j-2}}+\dfrac{1}{a_{2j-1}}\right )\cdot  
\left ( -\dfrac{1}{a_{2j-1}}+\dfrac{1}{a_{2j}}\right )$
and hence $\;|\delta _{j-1}|\leq \left |\dfrac{1}{a_{2j-2}a_{2j-1}} \right |
+\left |\dfrac{1}{a_{2j-1}a_{2j}} \right |+ \dfrac{1}{10}
\left | \dfrac{-1}{a_{2j-2}} +\dfrac{1}{a_{2j-1}} \right |\cdot
\left | \dfrac{-1}{a_{2j-1}} +\dfrac{1}{a_{2j}} \right |
\leq 1+1+\dfrac{1}{10}\cdot 2 \cdot 2=2.4.$
Also $|\delta _j | \leq 2.4 \;$ and thus $\gamma _j >|\delta _{j-1}|+
|\delta _j |.$
This proves the right inequality.
\end{proof}

\begin{rem1} $ 6 $ is the best integer upper bound.
For the proof see Remark 2 in section 8.
\end{rem1}

\section{Theorem 2: The case of real zeros}

\begin{thm2}
If $a_ja_{j+ 1}<0,\;$ then all zeros of $\;\Delta _K(t)\;$ are real and positive.
\end{thm2}
\begin{proof} 
We show that $\Delta _K(t) $ has a symmetric companion matrix.

\noindent
 Let $\;r=[2a_1,-2a_2,2a_3, \ldots, (-1)^{m-1}2a_m],\;$
where $ a_j>0. \;$ Then the Seifert matrix $ U $ is of the form 
\begin{equation*}
U= \left [
\begin{array}{ccccccc}
a_1 & 0 & & & & \\
-1& -a_2 & 1& & & \\
& 0 & a_3 & 0 &  & \\
& & -1 & -a_4 & 1 & \\
&&&&a_5&\\

& & & & \ddots &  
\end{array}
\right ]
\end{equation*}
Now 
\begin{equation*}
Ut-U^T= \left [
\begin{array}{ccccc}
a_1(t-1)&1&&&\\
-t&-a_2(t-1)& t & &\\
& -1 &a_3(t-1) & 1 &\\
&&& \ddots&
\end{array}
\right ]
\end{equation*}
We apply a series of transformations that don't change
the zeros of the determinant of the matrix.
First, multiply $ -1 $ on all even rows to get
\begin{equation*}
 \left [
\begin{array}{ccccc}
a_1(t-1)&1&&&\\
t&a_2(t-1)& -t & &\\
& -1 &a_3(t-1) & 1 &\\
&&& \ddots&
\end{array}
\right ]
\end{equation*}
Then multiply $ \dfrac{1}{\sqrt{a_1}} $ on the 1-st row and column,
 $ \dfrac{1}{\sqrt{a_2}} $ on the 2-nd row and column, and so on, to get
\begin{equation*}
M=\left [
\begin{array}{ccccc}
t-1 & \dfrac{1}{\sqrt {a_1a_2}} &&&\\
\dfrac{t}{\sqrt{a_1a_2}} & t-1 & \dfrac{-t}{\sqrt{a_2a_3}} &&\\
&\dfrac{-1}{\sqrt{a_2a_3}} & t-1 & \dfrac{1}{\sqrt{a_3a_4}} &\\
&&&\ddots
\end{array}
\right ]
\end{equation*}
with $\det(M)=\dfrac{1}{a_1a_2\cdots a_m} \det (Ut-U^t) .$
Now eliminate $ t $ from the off-diagonal line as follows:
multiply $-\dfrac{1}{\sqrt{a_1a_2}} $ on the 1-st row and add it 
to the 2-nd row,
multiply $ \dfrac{1}{\sqrt{a_2a_3}} $ on the 3-rd row and add it 
to the second row,
multiply $-\dfrac{1}{\sqrt{a_3a_4}} $ on the 3-rd row and add it 
to the 5-th row, etc,
i.e. multiply $ M $ by matrix $ P $ from the left:
\begin{equation*}
P= \left [
\begin{array}{ccccccc}
1 & 0 &&&&&\\
-\dfrac{1}{\sqrt{a_1a_2}}& 1&\dfrac{1}{\sqrt{a_2a_3}}&&&&\\
&0&1&0&&&\\
& & -\dfrac{1}{\sqrt{a_3a_4}}&1&\dfrac{1}{\sqrt{a_4a_5}}&&\\
&&&0&1&&\\
&&&&&& \ddots 
\end{array}
\right ]
\end{equation*}
Then $PM =$
\begin{equation*}
\left [
\begin{array}{ccccccc}
t-1 & \dfrac{1}{\sqrt{a_1a_2}}&0&0&&&\\
\dfrac{1}{\sqrt{a_1a_2}}&(t\!-\!1\!-\!\dfrac{1}{a_1a_2}\!-\!\dfrac{1}{a_2a_3})&
\dfrac{-1}{\sqrt{a_2a_3}}&\dfrac{1}{\sqrt{a_2a_3}\sqrt{a_3a_4}}&0&&\\
0&\dfrac{-1}{\sqrt{a_2a_3}}&t-1&\dfrac{1}{\sqrt{a_3a_4}}&0&0\\
0&\dfrac{1}{\sqrt{a_2a_3}\sqrt{a_3a_4}}&
\dfrac{1}{\sqrt{a_3a_4}}&(t\!-\!1\!-\!\dfrac{1}{a_3a_4}\!-\!\dfrac{1}{a_4a_5})&
\dfrac{-1}{\sqrt{a_4a_5}}&\dfrac{1}{\sqrt{a_4a_5}\sqrt{a_5a_6}}&\\
&0&0&\dfrac{-1}{\sqrt{a_4a_5}}&t-1&\dfrac{1}{\sqrt{a_5a_6}}&\\
&&&&&&\ddots 
\end{array}
\right ]
\end{equation*}
Since $PM=tE-A ,\; A \;$ is a companion matrix of $\Delta _K(t)$ and it is symmetric.
 So all its eigenvalues are real, and hence positive.
\end{proof}

\section{ Theorem 3: The case $a_ia_{i+1}\neq 1$} 
\begin{thm3}
Let $r=[2\eps _1a_1,2\eps _2a_2,\ldots ,2\eps _ma_m], \; $ where
$\; a_i>0,\;\eps _i=\pm1.\;$ If we don't have $\;a_i=a_{i+1}=1 $ and
 $\eps _i=\eps _{i+1}, $ then the zeros of $\Delta _{K(r)} $
satisfy inequality:
\begin{equation*}
-1< \Re (\alpha ).
\end{equation*} 
If, moreover, $a_j>1 $ for all $j$, then $\Re(\alpha )<3 $.
\end{thm3}
\begin{proof}
We find a positive definite (symmetric) matrix $ V $ such that 
$ V(E+A)+(E+A^T)V=W $ is positive definite. Let $ V $ be a diagonal matrix
 with elements $ a_1,a_2,\ldots ,a_m. $ Then multiplying $ (E+A) $
by $ V $ from the left is multiplying the $ i-$th row of $ (E+A) $
by $ a_i,\;i=1,\ldots ,m. $ Let $ l=[\frac{m}{2}]. $
Define $\eps _{ij}=\eps _i \eps _j$.

By \eqref{U^-1U^T} we have $ E+A= $ 
\begin{equation*}\label{E+A}
= \left [
\begin{array}{clclcccc}
2\hspace{2mm} & -\frac{\eps _1}{a_1} & 0 &&\ldots &&&\\
\vspace{2mm}
\frac{\eps _2}{a_2} \hspace{2mm}& 2\!-\!\frac{\eps _{12}}{a_1a_2}\!-\!
\frac{\eps _{23}}{a_2a_3}&\hspace{2mm}
-\frac{\eps _2}{a_2}\hspace{2mm} & \frac{\eps _{23}}{a_2a_3}& 0 & \ldots & &\\
\vspace{2mm}
0&\frac{\eps _3}{a_3}& 2 & -\frac{\eps _3}{a_3}& 0 & 0 & 0 & \ldots \\
\vspace{2mm}
 0 & \frac{\eps _{34}}{a_3a_4} & \frac{\eps _4}{a_4} & 2\!-
\!\frac{\eps _{34}}{a_3a_4}
\!-\!\frac{\eps _{45}}{a_4a_5}&\hspace{2mm} -\frac{\eps _4}{a_4} &
 \frac{\eps _{45}}{a_4a_5}& 0 & \ldots  \\
\vspace{2mm}
0& 0 &0& \frac{\eps _5}{a_5}& 2 & -\frac{\eps _5}{a_5}& 0& \ldots\\
\vspace{4mm}
&&& \ldots &&&\\
\end{array}
\right ]
\end{equation*}
where the last row is $ (0,\ldots, 0,\frac{\eps _m}{a_m},2) $ if $ m $ is odd,

\noindent
and $ (0,\ldots, 0, \frac{\eps _{m-1,m}}{a_{m-1}a_{m}},\frac{\eps _m}{a_m},
2\!-\!\frac{\eps _{m-1,m}}{a_{m-1}a_m},) $ if $ m $ is even.

\noindent
 Therefore $ V(E+A)=$
\begin{equation*}\label{V(E+A)}
= \left [
\begin{array}{clclcccc}
2a_1\hspace{2mm} & -\eps _1 & 0 &&\ldots &&&\\
\vspace{2mm}
\eps _2 \hspace{2mm}& 2a_2\!-\!\frac{\eps _{12}}{a_1}\!-\!\frac{\eps _{23}}{a_3}&\hspace{2mm}
-\eps _2\hspace{2mm} & \frac{\eps _{23}}{a_3}& 0 & \ldots & &\\
\vspace{2mm}
0&\eps _3& 2a_3 & -\eps _3& 0 & 0 & 0 & \ldots \\
\vspace{2mm}
 0 & \frac{\eps _{34}}{a_3} & \eps _4 & 2a_4\!-
\!\frac{\eps _{34}}{a_3}
\!-\!\frac{\eps _{45}}{a_5}&\hspace{2mm} -\eps _4 & \frac{\eps _{45}}{a_5}& 0 & \ldots  \\
\vspace{2mm}
0& 0 &0& \eps _5& 2a_5 & -\eps _5& 0& \ldots\\
\vspace{4mm}
&&& \ldots &&&\\
\end{array}
\right ]
\end{equation*}
where the last row is $ (0,\ldots,0,\eps _m,2a_m) \;$ if $ m $ is odd, and

\noindent
$ (0,\ldots,0,\eps _m,2a_m-\frac{\eps _{m-1,m}}{a_{m-1}})\; $ 
if $ m $ is even.
\vspace{2mm}
\noindent
 Further $ W=V(E+A)+(E+A^T)V=$
\begin{equation*}\label{W}
 \left [
\begin{array}{clclclcc}
4a_1 & -\eps _1+\eps _2 & 0 &&\ldots &&\\
\vspace{2mm}
-\eps _1+\eps _2 & 4a_2\!-\!\frac{2\eps _{12}}{a_1}\!-\!\frac{2\eps _{23}}{a_3}&
-\eps _2+\eps _3 & \frac{\eps _{23}+\eps _{34}}{a_3}&  \ldots & &\\
\vspace{2mm}
0&-\eps _2+\eps _3& 4a_3 & -\eps _3+\eps _4& 0 & 0 &  \ldots \\
\vspace{2mm}
 0 & \frac{\eps _{23}+\eps _{34}}{a_3} &-\eps _3+ \eps _4 & 4a_4\!-
\!\frac{2\eps _{34}}{a_3}
\!-\!\frac{2\eps _{45}}{a_5}& -\eps _4+\eps _5 & \frac{\eps _{45}+\eps _{56}}{a_5}
&  0 \\
\vspace{2mm}
0& 0 &0& -\eps _4+\eps _5& 4a_5 & -\eps _5+\eps _6& 0\\
\vspace{4mm}
0&0&0&\frac{\eps _{45}+\eps _{56}}{a_5}&-\eps _5+\eps _6&
4a_6-\frac{2\eps _{56}}{a_5}-\frac{2\eps _{67}}{a_7}&\ldots\\
&&& \ldots &&&\ddots\\
\end{array}
\right ]
\end{equation*}
The last row of $W$ is $ (\;0,\ldots, 0,-\eps _{m-1}+\eps _m,\:4a_m\;) $
if $ m $ is odd, and 

\noindent
$ (\;0,\ldots,0,\frac{\displaystyle \eps _{m-2,m-1}
+\eps _{m-1,m}}{\displaystyle a_{m-1}},\;-\eps _{m-1}+\eps _m,\:
4a_m-\frac{\displaystyle 2\eps _{m-1,m}}{\displaystyle a_{m-1}}\;) $
if $ m $ is even.

We eliminate the elements $-\eps _i+\eps _{i+1}$ : if $i$ is odd, multiply the 
$i $-th row  by $(\eps _i-\eps _{i+1})/4a_i $ and add to 
the $(i+1)$-th row. If $ i $ is even, multiply the $(i+1)$-th row by $(\eps _i-\eps _{i+1})/4a_{i+1}$ and add to the $i$-th row.
Similarly for columns.
In other words we consider the matrix $PWP^T,$  where 
\begin{equation*}\label{P}
P=\left [
\begin{array}{cccccc}
1&&&&&\\
\vspace{1mm}
\frac{\eps _1-\eps _2}{4a_1}&1&\frac{\eps _2-\eps _3}{4a_3}&&&\\
\vspace{1mm}
&&1&&&\\
\vspace{1mm}
&&\frac{\eps _3-\eps _4}{4a_3}&1&\frac{\eps _4-\eps _5}{4a_5}&\\
\vspace{1mm}
&&&&1&\\
&&&&&\ddots
\end{array}
\right ]
\end{equation*}
We have
\begin{equation*}\label{PWP^T} 
 PWP^T =\left [
\begin{array}{ccccccc}
4a_1&0&0&&&&\\
0&\alpha _2&0&\beta _2&0&0&\\
0&0&4a_3&0&0&\\
0&\beta _2&0&\alpha _4&0&\beta _4&\ldots \\
&&&0&4a_5&0&\ldots \\
&&&\beta _4&0&\alpha _6&\ldots\\
&&&&&&\ldots
\end{array}
\right ]\sim
\end{equation*}
\begin{equation}\label{oplus}
\sim \left [ 
\begin{array}{ccccc}
4a_1 &&&&\\
&4a_3 &&&\\
&&\ddots &&\\
&&&&4a_{2l\pm 1}
\end{array}
\right ] \oplus \left [
\begin{array}{ccclc}
\alpha _2&\beta _2&&&\\
\beta _2&\alpha_4&\beta _4&0&\\
&\beta _4&\alpha _6& \beta _6&\ldots \\
&&&\ddots &\\
&&&\beta _{2l-2}& \alpha _{2l}
\end{array} 
\right ]
\end{equation}
Here $l=[\frac{\displaystyle m}{\displaystyle 2}]$, and for $i=1,\ldots , l-1 $
\begin{equation*}
\alpha _{2i}=4a_{2i}-\frac{\displaystyle 2\eps _{2i-1,2i}}{\displaystyle a_{2i-1}}-
\frac{\displaystyle 2\eps _{2i,2i+1}}{\displaystyle a_{2i+1}}-\frac{\displaystyle (\eps _{2i-1}-\eps _{2i})^2}
{\displaystyle 4a_{2i-1}}-
\frac{\displaystyle (\eps _{2i}-\eps _{2i+1})^2}{\displaystyle 4a_{2i+1}}
\end{equation*}
\begin{equation}\label {alpha'}
=4a_{2i}-\frac{\displaystyle 3}{\displaystyle 2}\frac{\displaystyle \eps _{2i-1,2i}}{\displaystyle a_{2i-1}}-
\frac{\displaystyle 3}{\displaystyle 2}\frac{\displaystyle \eps _{2i,2i+1}}{\displaystyle a_{2i+1}}-\frac{\displaystyle 1}{\displaystyle 2a_{2i-1}}-
\frac{\displaystyle 1}{\displaystyle 2a_{2i+1}}, 
\end{equation}
\begin{equation*}
\beta _{2i}=\frac{\displaystyle \eps _{2i,2i+1}+\eps _{2i+1,2i+2}}{\displaystyle a_{2i+1}}-
\frac{\displaystyle (\eps _{2i}-\eps _{2i+1})(\eps _{2i+1}-\eps _{2i+2})}
{\displaystyle 4a_{2i+1}}
\end{equation*}
\begin{equation*}
=\frac{\displaystyle 3}{\displaystyle 4} 
\frac{\displaystyle (\eps _{2i,2i+1}+\eps _{2i+1,2i+2})}
{\displaystyle a_{2i+1}}
+\frac{\displaystyle \eps _{2i,2i+2}+1}{\displaystyle 4a_{2i+1}},
\end{equation*}
\begin{equation*}
\alpha _{2l}=4a_{2l}-\dfrac{3}{2}\dfrac{\eps _{2l-1,2l}}{a_{2l-1}}-
\frac{\displaystyle 3}{\displaystyle 2}\frac{\displaystyle \eps _{2l,2l+1}}{\displaystyle a_{2l+1}}-\frac{\displaystyle 1}{\displaystyle 2a_{2l-1}}-
\frac{\displaystyle 1}{\displaystyle 2a_{2l+1}},\;\; 
\text{ if } m \text{ is odd, }
\end{equation*} 
\begin{equation*}
\alpha _{2l}=4a_{2l}-\frac{\displaystyle 3}{\displaystyle 2}\frac{\displaystyle \eps _{2l-1,2l}}{\displaystyle a_{2l-1}}-
-\frac{\displaystyle 1}{\displaystyle 2a_{2l-1}},\;\; \text{ if } m  
\text{ is even.}
\end{equation*}

\noindent
Since all $ a_j \geq 1, $  it is not difficult to check
that if among $\;\frac{\displaystyle \eps _{2i-1}}{\displaystyle a_{2i-1}},
\;\frac{\displaystyle \eps _{2i}}{\displaystyle a_{2i}},\;\;
\frac{\displaystyle \eps _{2i+1}}{\displaystyle a_{2i+1}},\;\;$
 there are no two consecutive $ 1 $ or $ -1 $,
then the conditions of Positivity Lemma are satisfied: 
$\;(i) \;\alpha _{2i} > 0,\;$ 
$(ii) \;\alpha _{2i} \geq |\beta _{2i-2}|+|\beta _{2i}|,\; i=2, \ldots, l-1,\; $ and 
$\;(iii)\;\alpha _2> |\beta _2| $. 
If $ \beta _{2j}=0 $ then $\alpha _{2j+2}>|\beta _{2j+2}|$.
So the second matrix in \eqref{oplus} is positive definite and 
so is $ W $. 

The proof of inequality $ \Re(\alpha )<3 $ in the case $a_j>1$ for $j=
1,2,\ldots,m $,  is similar to the proof 
of Theorem 1 for $q=3$.
\end{proof}

\section{Theorem 4: The case of fibered knots}

Consider a fibered two-bridge knot $K(r)$ 
with $r=[2a_1,2a_2,\ldots ,2a_m]$
\begin{equation*}
\begin{array}{ccccc}
=[&\underbrace{ 2, \ldots , 2,} &  \underbrace{-2, \ldots, -2},&
\ldots,&
\underbrace{(-1)^{m-1}2,\ldots , (-1)^{m-1}2}\;] \\
&k_1 &  k_2&& k_m
\end{array}.
\end{equation*}

\noindent
The following theorem is a corollary of Theorem 3:
\begin{thm4}
If $k_j=1  \text{ or } 2, \; j=1, \ldots m, \; $   
then $ -1 < \Re (\alpha )  $.
\end{thm4}
\begin{proof}
At least one of $\eps _{2i-1,2i}, \;\eps _{2i,2i+1} $ in \eqref{alpha'} is negative.
So by \eqref{alpha'}
\begin{equation*}
 \alpha _{2i} = 
\begin{cases}
 3 & \text{ if } \eps _{2i-1}\neq \eps _{2i+1}\\
 6 & \text{ if } \eps _{2i-1}=\eps _{2i+1}\neq \eps _{2i}
\end{cases}
\end{equation*}
While 
\begin{equation*}
\beta _{2i}= 
\begin{cases} 
0 &\text{ if }\eps _{2i}\neq \eps_{2i+2}\\
-1 &\text{ if } \eps _{2i}=\eps _{2i+2}\neq \eps _{2i+1}
\end{cases}
\end{equation*}
and similarly $\beta _{2i-2} =0 $ or $-1.$
So the conditions of Positivity Lemma are satisfied, which proves the 
inequality. 
\end{proof}  

\section{Theorem 5: The case $\;a_i=\pm c\;$}
\begin{thm5}
Let $\;r_m=[2c,-2c,\ldots ,(-1)^{m-1}2c],\;\;c>0,\;m\geq 1.\;$
Then all zeros of $\;\Delta _{K(r_m)}\; $ satisfy inequality:
\begin{equation*}
(\dfrac{\sqrt{1+c^2}-1}{c})^2<\alpha <(\dfrac{\sqrt{1+c^2}+1}{c})^2.
\end{equation*}
\end{thm5}
\begin{proof}  
By \eqref{U} a  Seifert matrix for $\; K(r_m) \;$ is 
\begin{equation*}
U=\left [
\begin{array}{cccccc}
c&0&&&&\\
-1&-c&1&&&\\
&0&c&0&&\\
&&-1&-c&1&\\
&&&&\ddots &
\end{array}
\right ]
\end{equation*}

\noindent
Let $ P_0(t)=1,\; P_1(t)=c(t-1),\;
P_m(t)=(-1)^{[\frac{m}{2}]} \det (tU-U^T)=$
\begin{equation*}
=(-1)^{[\frac{m}{2}]}  \det \left [
\begin{array}{cccccc}
c(t-1)&1&&&&\\
-t&c(-t+1)&t&&&\\
&-1&c(t-1)&1&&\\
&&-t&c(-t+1)&t&\\
&&&&\ddots&
\end{array}
\right ]\;=\;
\end{equation*}
\begin{equation*}  
= \det \left [ 
\begin{array}{ccccc}
c(t-1)&1&&&\\
t&c(t-1)&-t&&\\
&-1&c(t-1)&1&\\
&&t&c(t-1)&\\
&&&\ddots &
\end{array}
\right ]
\end{equation*}
Then $P_m(t)=\pm \Delta _{K(r_m)}(t)$, and  $P_m(t)\;$ satisfy a recurrence equation:
\begin{equation}\label{c1P_m}
P_m(t)=c(t-1)P_{m-1}(t)-tP_{m-2}(t), \; m\geq 2.
\end{equation}
Since $K(r_{2m+1})$  is a 2-component link, we can
write $P_{2m+1}(t)=(t-1)Q_{2m}(t)$.
Note $Q_0(t)=c$. Then from \eqref{c1P_m} we have
\begin{equation}\label{c2P_2m}
P_{2m}(t)=c(t-1)^2Q_{2m-2}(t)-tP_{2m-2}(t).
\end{equation}
Also,  $$P_{2m+1}(t)=c(t-1)P_{2m}(t)-tP_{2m-1}(t) \Longrightarrow $$
$$(t-1)Q_{2m}(t)=c(t-1)P_{2m}(t) - t(t-1)Q_{2m-2}(t)\Longrightarrow $$
\begin{equation}\label{c3Q_2m}
Q_{2m}(t)=cP_{2m}(t)-tQ_{2m-2}(t).
\end{equation}
\vspace{2mm}
 Then \eqref{c2P_2m}  
and \eqref{c3Q_2m} imply
$$t^{-m}P_{2m}(t)=t^{-m}c(t-1)^2Q_{2m-2}(t)-t^{-(m-1)}P_{2m-2}(t)$$ 
and 
$$t^{-m}Q_{2m}(t)=ct^{-m}P_{2m}(t)-t^{-(m-1)}Q_{2m-2}(t).$$ 
Let $x=t+\dfrac{ 1}{ t},\;$ and write 
$\;\phi _m(x)=t^{-m}P_{2m}(t), \;\;\psi _m(x)=t^{-m}Q_{2m}(t). $ Then
\begin{equation}\label{c4phi_m}
\phi _m(x)=c(x-2)\psi _{m-1}(x)-\phi _{m-1}(x),
\end{equation}
\begin{equation}\label{c5psi_m}
\psi _m(x)=c\phi _m(x)-\psi _{m-1}(x).
\end{equation}
Note $\phi _0(x)=1, \:\psi _0(x)=c. \;$
Since \eqref{c4phi_m} $\Longrightarrow 
c(x-2)\psi_{m-1}(x)=\phi_m(x)+\phi_{m-1}(x),$ 
from \eqref{c5psi_m} we see:
$$c(x-2)\psi _m(x)=c^2(x-2)\phi _m(x)-c(x-2)\psi _{m-1}(x)\Longrightarrow  $$
$$\phi _{m+1}(x)+\phi _m(x)=c^2(x-2)\phi _m(x)-
(\phi _m(x)+\phi _{m-1}(x)) \Longrightarrow \;$$
\begin{equation}\label{c6phi_m+1}
\phi _{m+1}(x)=(c^2x-(2c^2+2))\phi _m(x)-\phi _{m-1}(x)
\end{equation}
Similarly, using \eqref{c4phi_m} and \eqref{c5psi_m}, we have
\begin{equation}\label{c7}
\psi _m(x)=(c^2x-(2c^2+2))\psi _{m-1}(x)-\psi _{m-2}(x)
\end{equation}
Let $\;y=c^2x-(2c^2+2). \;$
Write $\;\phi _m(x)=\lambda _m(y) \;$ and $\;\psi _m(x)=\mu _m(y).\;$
Then from \eqref{c6phi_m+1} and \eqref{c7} we have, for $m \geq 2$,
\begin{equation*}\label{c8}
\lambda _m(y)=y\lambda _{m-1}(y)-\lambda _{m-2}(y)
\end{equation*}
\begin{equation*}\label{c9}
\mu _m(y)=y\mu_{m-1}(y)-\mu _{m-2}(y),
\end{equation*}
where $ \lambda _0=1,\;\lambda _1=y+1,\;\lambda _2=y^2+y-1,\;
\mu _0=c,\;\mu _1=cy,\;\mu _2=c(y^2-1).\;$
It is easy to see that for $ m \geq 1,\; 
\lambda _m=\dfrac{1}{c}(\mu _m+\mu_{m-1}).$
Now  let $\;f_m(y)\;$ be a Fibonacci polynomial defined in \cite{K}:
$f_1(y)=1,\; f_2(y)=y \;$ and for $ m \geq 3, $
\begin{equation*}\label{c9f_m}
f_m(y)=yf_{m-1}(y)+f_{m-2}(y).
\end{equation*}
Then we can show by induction that for $m \geq 0$,
\begin{equation*}
i^{-m}f_{m+1}(iy)=\dfrac{1}{c}\mu _m(y).
\end{equation*}
It is known (see \cite{K}, p.477) that the 
zeros of $f_{m+1}(y)$ are $y_k=2i\cos\frac{k\pi }{m+1}\:,\;$

\noindent
$k=1,2,\ldots,m.$
Therefore, the zeros of $\; \mu _m(y) \; $ are
\begin{equation*}\label{c11y_k}
y_k^{(m)} = 2\cos \frac{k\pi }{m+1},\;\;k=1,2,\ldots ,m.
\end{equation*}
Next we look at the zeros of $\;\lambda _m(y).\;$ Since
$y_k^{(m-1)}=2\cos \frac{k\pi }{m},\;\;k=1,2,\ldots ,m-1,$ are all the zeros
of $\;\mu _{m-1}(y),\;$ and for any $\;k\;$ 
\begin{equation}\label{c12y}
 y^{(m-1)}_{k+1} < y_{k+1}^{(m)} < y_k^{(m-1)}< y_k^{(m)},
\end{equation}
 there exists exactly one zero of $\;\mu _m(y)\;$ between
neighboring two zeros of $\;\mu _{m-1}(y),\;$ and also there exists
exactly one zero of $\;\mu _{m-1}(y) \;$ between neighboring two zeros
of $\; \mu _m(y) \;$ (see Fig.3). By induction we check that
\begin{equation}\label{c13mu_2m}
\mu _{2m}(-2)=(2m+1)c \;\;\text{ and }\;\;\mu _{2m+1}(-2) =-(2m+2)c.
\end{equation}

\noindent
Now, the zeros of $\lambda _m(y) $ occur at the intersections of two curves 
 $c_1:\;z=(-1)^m\mu _m(y)  \;$ and $\;
c_2:\;z= (-1)^{m-1}\mu _{m-1}(y).\;$
By \eqref{c12y} there are $\;m-1 \;$ zeros in $\;(y^{(m)}_{m-1},2), \;$ and by
\eqref{c13mu_2m} two curves intersect in $\;(-2,y_{m-1}^{m}).\;$
Therefore there are
exactly $\;m\;$ real zeros in $\;(-2,2).\;$
Since $y=c^2x-(2c^2+2),\;x=\dfrac{y+(2c^2+2)}{c^2}\;$ and the zeros of 
$\;\phi _m(x)\;$ and $\;\psi _m(x)\;$ are in the interval 
$\;(2,\;2+\dfrac{4}{c^2}),$ and hence all zeros of $\;P_{2m}(t)\;$ and 
$\;Q_{2m}(t)\;$ satisfy inequality:
\begin{equation*}
\dfrac{1}{q}=(\dfrac{\sqrt{1+c^2}-1}{c})^2< \alpha <q=(\dfrac{\sqrt{1+c^2}+1}
{c})^2.\qedhere
\end{equation*}
\end{proof}

\begin{figure}[h]
\includegraphics[height=10cm]{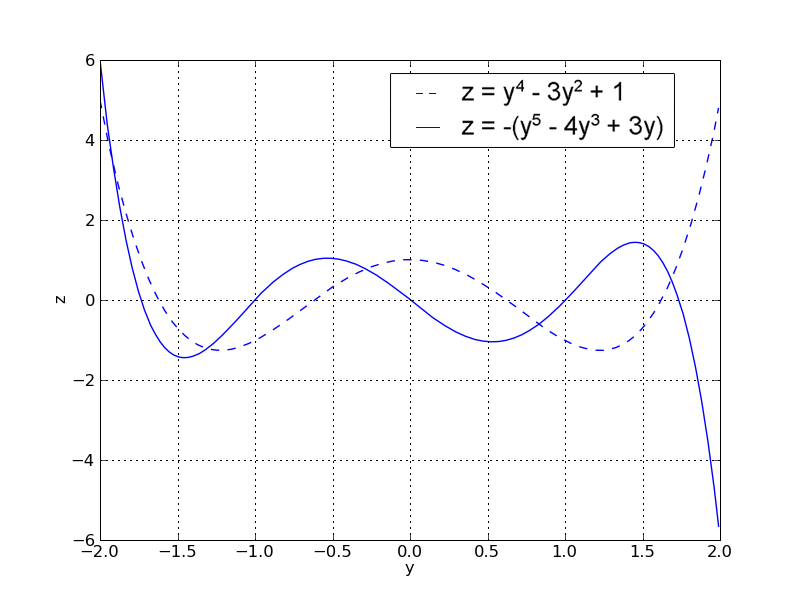}

Figure 3.
\label{graph}
\end{figure}
\begin{cor1}
If $\;c \rightarrow \infty,\;$ then the zeros 
of $P_{2m}(t),\; Q_{2m}(t),$ which are the zeros of Alexander polynomials, tend to $1$.
\end{cor1}
\begin{rem2}
For $c=1 $ and large enough $m$ we can find a zero $\alpha $ of 
$P_{2m}(t)$
arbitrarily close to $q=3+\sqrt{8}$. It is quite likely that $3+\sqrt{8}$
is the upper bound of the real part of the zeros.
\end{rem2}
\begin{proof}
Since the zeros of $(-1)^{m-1} \mu _{m-1}(y) \;$ and $\;(-1)^m \mu _m(y)\;$ satisfy inequality $y_2^{(m)} < y_1^{(m-1)} < y_1^{(m)},$ there is a zero of $\lambda _m$
greater than $y_2^{(m)}$, where $y_2^{(m)}=2\cos \dfrac{2\pi }{m+1}$.
So there is a zero of $\phi _m(x)$ arbitrarily close to $6$, hence a zero of $P_{2m}(t)$ arbitrarily close to $3+\sqrt{8}$.
\end{proof}

\section{Open questions} 

Let us finish with several open questions:
\smallskip

\noindent 1) Is there an upper bound of the real part of zeros of
the Alexander polynomials of general alternating knots ?
Recently Hirasawa observed(2010) that each of the following alternating
12 crossing knots $12a_{0125}$ and $12a_{1124}$ has a real zero, 6.90407... and
7.69853... respectively. Therefore an upper bound, if exists, is larger
than 7.
\medskip

\noindent 2) Given m, does there exist an upper bound q(m) of the real part
of zeros of the Alexander polynomials of degree m of alternating knots ?
\medskip

\noindent 3) Is there a version of Conjecture 1 for 
non-alternating knots ? 

\noindent
Notice that  Conjecture 1 does not hold for homogeneous knots (defined  
in  [Cr]).
Hirasawa showed (2010) that a
non-alternating knot $10_{152}$ is a closure of a positive 3-braid and hence
it is a homogeneous knot, but the Alexander polynomial has a real zero
$\alpha=-1.85... $ 
\medskip

\noindent 
4) Characterize alternating knots whose zeros of the Alexander polynomial are real. In particular, is the converse of Theorem 2 true
for one component two-bridge knots?
\bigskip

{\bf Aknowledgements} :
we are grateful to Misha Lyubich for a helpful reference and to 
Yun Tao Bai for his help with figures.

\end{document}